\documentclass[preprint]{elsarticle}

\usepackage{lineno,hyperref}
\journal{Discrete Applied Mathematics}

\usepackage{enumerate}
\usepackage{amsmath,amssymb,amsthm}
\usepackage{graphicx}
\usepackage{tikz}
\usepackage{color}

\newtheorem{thm}{Theorem}[section]

\newtheorem{prop}[thm]{Proposition}

\newtheorem{defi}[thm]{Definition}

\newdefinition{rmk}{Remark}
\newdefinition{example}{Example}

\newenvironment{pf}{\begin{proof}}{\end{proof}}

\begin{document}
%%%%%%%%%%%%%%%%%%%%%%%%%%%%%%%%%%%%%%%%%%%%%%%%%%%%%%%%%%%%%%%%%%%%%
%%%%%%%%%%%%%%%%%%%%%%%%%%%%%%%%%%%%%%%%%%%%%%%%%%%%%%%%%%%%%%%%%%%%%

\begin{frontmatter}

\title{Uniform forcing and immune sets in graphs and hypergraphs\tnoteref{t1}}
\tnotetext[t1]{Partially supported by the Ministerio de Ciencia e Innovaci\'on\mbox{$/$}Agencia Estatal de Investigaci\'on, Spain, and the European Regional Development Fund under project  PGC2018-095471-B-I00; and by AGAUR from the Catalan Government under project 2017SGR--1087.}

\author[1]{Josep F\`abrega\corref{cor1}}
\ead{josep.fabrega@upc.edu}
\author[1]{Jaume Mart\'{\i}-Farr\'e}
\ead{jaume.marti@upc.edu}
\author[1]{Xavier Mu\~{n}oz}
\ead{xavier.munoz@upc.edu}

\address{Departament de Matem\`atiques, Universitat Polit\`ecnica de Catalunya, Jordi Girona 1--3, 08034, Barcelona, Spain}

%\affiliation[1]{organization={Departament de Matem\`atiques, Universitat Polit\`ecnica de Catalunya}, addressline={Jordi Girona 1--3}, postcode={08034}, city={Barcelona}, country={Spain}}

\cortext[cor1]{Corresponding author}

\begin{abstract}
Zero forcing is an iterative coloring process on a graph that has been widely used in such  different areas as the modelling of propagation phenomena in networks and the study of minimum rank problems in matrices and graphs. This paper deals with zero forcing on hypergraphs. (Representing a network by a hypergraph allows us to account for its community structure and for  more general interactions between different subsets of nodes.) 

We consider two natural generalizations to hypergraphs of zero forcing on graphs (one of them already known) and, for each one of these generalizations  we look into two clutters that play a significant role in the forcing  process: the clutter of minimal forcing sets and the one of minimal immune sets. A formulation of immune sets in terms of neighbourhoods (hence without making reference to the iterative zero forcing process) is presented, highlighting the different behaviour of the distinct forcing rules. Moreover, we obtain the families of minimal forcing and minimal immune sets in the case of complete hypergraphs and we provide a full characterization of forcing and immune uniform clutters, both in the graph and in the hypergraph case.
\end{abstract}

\begin{keyword}
hypergraphs, zero forcing, propagation in networks.
\end{keyword}

\end{frontmatter}

%\linenumbers
%%%%%%%%%%%%%%%%%%%%%%%%%%%%%%%%%%%%%%%%%%%%%%%%%%%%%%%%%%%%%%%%%%%%%
\section{Introduction.}

The spread of disease in a population, the dissemination of information in a social network, or the cascading failure in a power electrical or computer network, are examples of propagation phenomena in networks which can be modelled by considering the following scenario: the network is represented by a simple undirected graph and the propagation phenomena by means of a discrete time dynamic forcing process acting on the vertices of the graph.  At any given time each vertex is in one of two possible states (black and white). At the beginning of the process, some vertices of the graph are colored black, while the rest of the vertices are initially colored white. At each step, one white vertex can become black by applying a well-defined forcing rule; and this forcing procedure is iteratively applied until no more changes of color are possible. See, for instance, the book \cite{BaBa-2008} for a  comprehensive overview of dynamical processes on networks.

In this paper we consider a more general framework in which the network is modelled by a hypergraph, hence allowing to account both for the community structure of the network and for more general interactions between different subsets of nodes than those taken into consideration by the adjacency relation in a graph. A \emph{hypergraph} $\mathcal{H}$ on a finite set $\Omega$ is a pair $(V(\mathcal{H}),\mathcal{E}(\mathcal{H}))$ where $V(\mathcal{H})=\Omega$ is the set of vertices of the hypergraph and $\mathcal{E}(\mathcal{H})=\{E_1,\ldots,E_m\}$ is the set of hyperedges $E_i\subseteq \Omega$, $1\leqslant i\leqslant m$, see for example the classic book \cite{Be-1973}. In addition, we assume that $\mathcal{E}(\mathcal{H})$ is a \emph{clutter}, that is,  none of the hyperedges  contains another (clutters are also known as \emph{antichains} or \emph{Sperner systems}, see for instance \cite{Ju-2011}). Observe that if $|E_i|=2$ for all $i$, $1\leqslant i\leqslant m$, then $\mathcal{H}$ is just a simple graph without isolated vertices. So, since any graph can be seen as a hypergraph, our model incorporates the graph representation of the network. 

Besides, in this work we consider different \emph{zero forcing rules} in order to model the propagation phenomena. In the graph case, zero forcing was introduced in \cite{AIM-2007, BuGi-2007} and works as follows: at each forcing step, a black vertex of a graph with exactly one white adjacent vertex will force this white vertex to become black. In the hypergraph case, the forcing rule will be such that at each forcing step a black subset of vertices of a certain hyperedge will force all the white vertices of this hyperedge to become black, see for instance \cite{BeRyFa-2018} for one possible definiton of zero forcing on hypergraphs. Again, the sequence of forcing steps takes place until no more color switches can happen. 

An initial set of black vertices such that all the vertices of the graph or the hypergraph become black after repeatedly applying the zero forcing rule is named a \emph{zero forcing set}. In the case of graphs, the cardinality of a smallest zero forcing set is known as the zero forcing number and this parameter has been heavily studied, see for instance the recent paper \cite{KaKaSu-2019} and its references. The zero forcing process has also been applied to minimum rank problems in matrices and graphs, see for example \cite{FaHo-2007} for a detailed discussion. In what follows, we will refer to zero forcing sets simply as forcing sets.

Although the hypergraph model of the network has not been so widely considered as the graph model, it has been previously taken into consideration. For instance, in \cite{BoKaSi-2016} the authors consider the classical SIS model of epidemic propagation in the case that the network is represented by a hypergraph; in \cite{Er-2019} the zero forcing generalization introduced in \cite{BeRyFa-2018} is accounted to improve the upper bound of the infection number of a hypergraph (which is the cardinality of a smallest forcing set); and in \cite{Ho-2020}, the value of the hypergraph zero forcing number and maximum nullity are determined for various families of uniform hypergraphs.

In our work we consider two natural generalizations of zero forcing on graphs. One of them is the definition considered in \cite{BeRyFa-2018} and the other one is introduced in our paper. For each one of these zero forcing rules on hypergraphs, we look into the clutters of \emph{minimal forcing sets} and \emph{minimal immune sets}. As we will see, these clutters play a remarkable role in the forcing process.  Roughly speaking, an immune set is a subset of white vertices whose color remains unchanged during the whole zero forcing process, given that initially the rest of the vertices of the hypergraph have been colored black.

Our main goals are, on the one hand, to provide general results on the clutters of minimal forcing and minimal immune sets for arbitrary hypergraphs (and hence also for graphs), and, on the other hand, to give more specific results in the particular case of the family of complete hypergraphs and uniform clutters. 

The paper is organized as follows. In Section~\ref{ZeroForcing} we consider the two different natural generalizations to hypergraphs of zero forcing on graphs (Definitions~\ref{rule R1} and \ref{rule R2}). For each zero forcing rule discussed in the paper we introduce the clutters of minimal forcing and minimal immune sets (Definitions~\ref{forcing} and \ref{immune}), and we relate them by using transversals (Proposition~\ref{transversal}). In Section~\ref{Neighbourhood characterization}, a nice characterization of the immune sets  formulated in terms of neighbourhoods is given, without reference to the iterative zero forcing process (Proposition~\ref{hyper-immune}). In the hypergraph case this neighbourhood characterization is different for each zero forcing rule, but both coincide in the graph case. Finally, in Section~\ref{clutters} we study the role that uniform clutters play on zero forcing on hypergraphs. So  we obtain on the one hand the families of minimal forcing and minimal immune sets in the case of complete hypergraphs (Propositions~\ref{uni-R2} and \ref{I-Uniform}), and on the other hand we provide a complete characterization of forcing and immune uniform clutters; that is, given  a uniform clutter we answer the question of the existence of a graph or a hypergraph such that its family of minimal forcing or minimal immune sets with respect to the considered zero forcing rule  is precisely the given uniform clutter (Theorems~\ref{s-uniformes-th-R1}, \ref{s-uniformes-th1} and \ref{s-uniformes-th2}).

 %%%%%%%%%%%%%%%%%%%%%%%%%%%%%%%%%%%%%%%%%%%%%%%%%%%%%%%%%%%%%%%%%%%%%%%%%%%
\section{Zero forcing rules on hypergraphs. Forcing and immune sets.}
\label{ZeroForcing}

Let $\mathcal{R}$ denote the forcing rule defined on the hypergraph $\mathcal{H}$ representing the network under consideration. We suppose that at each forcing step, and according to some well-defined condition formulated by $\mathcal{R}$, some subset of black vertices in certain hyperedge of $\mathcal{H}$ will force all the white vertices of this hyperedge to become black. The forcing rule $\mathcal{R}$ is iteratively applied until no more color changes are possible. In this paper we deal specifically with the \emph{zero forcing} rule defined for graphs as follows: 

\begin{defi}[Zero forcing rule $\mathcal{R}_0$, \cite{AIM-2007,BuGi-2007} ] 
\label{zero-forcing-graphs}
At each forcing step, a black vertex of the graph with exactly one white adjacent vertex will force this white vertex to become black. The zero forcing rule is then iteratively applied until no more changes of color are possible. 	
\end{defi}

An initial set of black vertices of the graph $G$ such that all the vertices of $G$ become black after repeatedly applying $\mathcal{R}_0$ is called a forcing set.  
%%%%%%%%%%%%%%%%%%%%%%%%%%%%%%%%%%%%%%%%%%%%%%%%%%%%%%%%%%%%%%%%%%%%%%%%%%%%%%%%%
\begin{figure}[htb]
\begin{center}
\begin{tikzpicture}[scale=0.55]
%\draw [help lines] (0,0) grid (8,8);
\draw[black,fill] (4,0) circle(4pt);
\draw[black,fill] (4,8) circle(4pt);
\draw[black,fill] (0,4) circle(4pt);
\draw[black,fill] (8,4) circle(4pt);
\draw[black,fill] (1.17,1.17) circle(4pt);
\draw[black,fill] (8-1.17,8-1.17) circle(4pt);
\draw[black,fill] (8-1.17,1.17) circle(4pt);
\draw[black,fill] (1.17,8-1.17) circle(4pt);
\draw[thick] (4,8)--(8-1.17,8-1.17)--(8,4)--(8-1.17,1.17)--(4,0)--(1.17,1.17)--(0,4)--(1.17,8-1.17)--(1.17,1.17)--(8-1.17,8-1.17)--(8-1.17,1.17);
 \draw[thick] (4,0)--(0,4);
\draw[thick] (4,8)--(8,4);
\draw[thick] (0,4)--(8,4);
\node[above] at (4,8.2) {$1$};
\node[right] at (8-1,8-1) {$2$};
\node[right] at (8.2,4) {$3$};
\node[right] at (7,1) {$4$};
\node[below] at (4,-0.2) {$5$};
\node[left] at (1,0.95) {$6$};
\node[left] at (-0.2,4) {$7$};
\node[left] at (1,7.05) {$8$};

\end{tikzpicture}
\end{center}
\caption{The sets $\{1,2,7\}$ and $\{3,5,6,7\}$ are forcing sets.}
\label{fig-def-zf}
\end{figure}
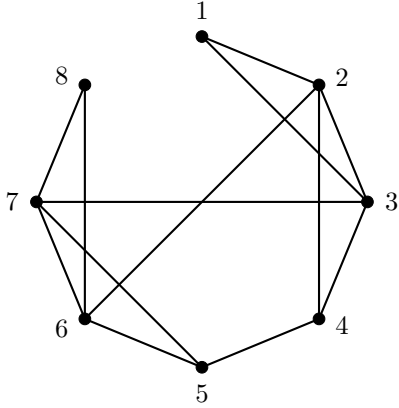
%%%%%%%%%%%%%%%%%%%%%%%%%%%%%%%%%%%%%%%%%%%%%%%%%%%%%%%%%%%%%%%%%%%%%%%%%%%%%%%%%
For example, in the graph $G$ shown in Figure~\ref{fig-def-zf}, the sets $\{1,2,7\}$ and $\{3,5,6,7\}$ are examples of forcing sets. 

For hypergraphs one can consider several natural generalizations of Definition~\ref{zero-forcing-graphs}. Let us remind first the notion of adjacency in hypergraphs. Given two subsets of vertices $X$ and $Y$ of a hypergraph $\mathcal{H}$ we say that $X$ and $Y$ are adjacent if there is a common hyperedge $E$ such that $X\subseteq E$ and $Y\subseteq E$. In particular, if $X$ and $\{v\}$ are adjacent for some vertex $v$, we say that $v$ is adjacent to the set $X$. 

A first generalization for hypergraphs, $\mathcal{R}_1$,  of zero forcing on graphs was introduced in \cite{BeRyFa-2018}. 

\begin{defi}[Forcing rule $\mathcal{R}_1$, \cite{BeRyFa-2018}]
\label{rule R1}
At each forcing step, a black subset $X$ of a hyperedge $E$, which is such that outside of $E$ there is not a white vertex adjacent to $X$, will force all the white vertices of $E$ to become black. This rule is iteratively applied until no more changes are possible.  
\end{defi}

Observe that if a graph $G$ is seen as a hypergraph $\mathcal{H}$, then the zero forcing rule $\mathcal{R}_1$ on $\mathcal{H}$ coincides with the zero forcing rule $\mathcal{R}_0$ on $G$.

In rule $\mathcal{R}_1$ the emphasis is put in the fact that the only white vertices adjacent to $X\subseteq E$ belong to the same hyperedge $E$, as happens in the graph case where the hyperedges are just the edges of the graph. But  in a graph we can adopt another perspective for zero forcing: a black vertex $b$ adjacent to a white vertex $w$ forces $w$ to become black if $b$ does not belong to any other edge with a white vertex. From this alternative point of view, let us  introduce a new natural generalization of zero forcing for hypergraphs, $\mathcal{R}_2$, which is different of the former one. As for $\mathcal{R}_1$, if a graph $G$ is seen as a hypergraph $\mathcal{H}$, then the zero forcing rule $\mathcal{R}_2$ on $\mathcal{H}$ coincides also with the zero forcing rule $\mathcal{R}_0$ on $G$. 

\begin{defi}[Forcing rule $\mathcal{R}_2$]
\label{rule R2}
At each forcing step, a black subset $X$ of a hyperedge $E$, which is such that $X$ is not contained in any other hyperedge with white vertices, will force all the white vertices of $E$ to become black. This rule is iteratively applied until no more changes are possible. 
\end{defi}

In accordance with the definitions of $\mathcal{R}_1$ and $\mathcal{R}_2$, it is clear that if some white vertex $v$ of a particular hyperedge can be forced to become black by applying the rule $\mathcal{R}_2$, then the color of $v$ can also be changed by $\mathcal{R}_1$. Nevertheless, the two rules are not equivalent. A simple example showing this fact is, for instance, the following one. Let $\Omega=\{1,2,3,4\}$ and consider the hypergraph $\mathcal{H}$ on $\Omega$ with hyperedges $\{1,2,3\}$, $\{1,2,4\}$ and $\{1,3,4\}$. Then we can easily check that if initially the set of black vertices is $\{1,2\}$, then at the end of the forcing process vertices $3$ and $4$ have become black by $\mathcal{R}_1$, whereas by rule $\mathcal{R}_2$ no change of color is possible.

Observe that for $i=0,1,2$, the zero forcing rule $\mathcal{R}_i$ satisfies the following property. Let $B\subseteq B'\subseteq V(\mathcal{H})$ be two subsets of black vertices, let $A\subseteq B$ and $A'\subseteq B'$ be two subsets such that $A\subseteq A'$, and let $E$ be a hyperedge such that $A'\subseteq E$. Therefore if according to $\mathcal{R}_i$ the white vertices of $E$ change its color  when we consider the black set $A$, then according to $\mathcal{R}_i$ these vertices also change its color when we consider the black set $A'$. Thus given an initial set $B$ of black vertices, if  $\mathcal{R}_i$ is iteratively applied until no more color changes are possible, the final set of black vertices resulting in $\mathcal{H}$ is unique and does not depend on the particular sequence of forcing steps performed on the hypergraph. Let us denote by ${\cal R}_i^{\ast}(B)$ this final black set. Clearly, $B\subseteq {\cal R}_i^{\ast}(B)$. Moreover, notice that if $B\subseteq B'$, then ${\cal R}_i^{\ast}(B) \subseteq {\cal R}_i^{\ast}(B')$. Let us refer to this last property by saying that the rule $\mathcal{R}_i$ is \emph{increasing}.

When the forcing process is finished, either all the vertices of the graph or the hypergraph have become black or there is still a non-empty remaining subset of white vertices. This observation motivates the following two definitions that we state in general for a hypergraph $\mathcal{H}$ on a finite set $\Omega$ and for the zero forcing rule $\mathcal{R}_i$, $i=0,1,2$.

\begin{defi}
\label{forcing}
Let $\mathcal{H}$ be a hypergraph on a finite set $\Omega$. An $\mathcal{R}_i$-forcing set of $\mathcal{H}$ is a non-empty subset $F\subseteq \Omega$ such that $\mathcal{R}^\ast_i(F)=\Omega$. Let us denote by  $\mathcal{F}_{i}(\mathcal{H})$ the family of the inclusion-minimal $\mathcal{R}_i$-forcing sets of $\mathcal{H}$; that is, $F\in \mathcal{F}_{i}(\mathcal{H})$ if and only if $F$ is an $\mathcal{R}_i$-forcing set and none of its proper subsets is an $\mathcal{R}_i$-forcing set. 
\end{defi}

\begin{defi}
\label{immune}
Let $\mathcal{H}$ be a hypergraph on a finite set $\Omega$. An $\mathcal{R}_i$-immune set of $\mathcal{H}$ is a non-empty subset $I\subseteq \Omega$ such that $\mathcal{R}^\ast_i(\Omega \setminus I)=\Omega \setminus I$. Let us denote by $\mathcal{I}_{i}(\mathcal{H})$ the family of the inclusion-minimal $\mathcal{R}_i$-immune sets of $\mathcal{H}$; that is, $I\in \mathcal{I}_{i}(\mathcal{H})$ if and only if $I$ is an $\mathcal{R}_i$-immune set and none of its proper subsets is an $\mathcal{R}_i$-immune set. 
\end{defi}

\begin{rmk}
In the case of zero forcing in graphs, the concept of immune set has been previously considered, in particular under the name of \emph{fort}. So, in \cite{FaHi-2018}, the authors define a fort $J$ of a graph as a non-empty subset of vertices such that no vertex outside of $J$ is adjacent to exactly one vertex of $J$. Thus, a fort is a set of vertices that do not get forced by some initial set of black vertices. It is worth noting that, in the case of graphs, the relationship between forts and immune sets will be more apparent in Section~\ref{Neighbourhood characterization} with the neighbourhood characterization of immune sets (see Remark~\ref{graph-immune} and Proposition~\ref{graph-immune-bis}).
\end{rmk}

If it is clear from the context which is the zero forcing rule $\mathcal{R}_i$ acting on $\mathcal{H}$, we can omit its reference and say plainly forcing set or immune set. 

\begin{rmk}
Since the zero forcing rule $\mathcal{R}_i$ is increasing, if $F$ is an $\mathcal{R}_i$-forcing set, then any superset $F'\supseteq F$ is also an $\mathcal{R}_i$-forcing set. This property is not true for immune sets. For instance, consider the zero forcing rule $\mathcal{R}_2$ in  the hypergraph $\mathcal{H}$ on $\Omega=\{1,2,3,4\}$ with hyperedges $\{1,2\}$, $\{1,3\}$, $\{2,3\}$ and $\{3,4\}$. We can easily check that $\{1,2\}$ is an immune set but $\{1,2,3\}$ is not. However, notice that if $I$ is an $\mathcal{R}_i$-immune set, then at the end of the forcing process all the vertices of $I$ remain white no matter which is the initial subset of black vertices $B\subseteq V(\mathcal{H})\setminus I$. 
\end{rmk}

By considering $\mathcal{F}_i(\mathcal{H})$ and $\mathcal{I}_i(\mathcal{H})$ some of our results can be formulated in terms of transversals (Proposition~\ref{transversal}). 

We use the following well-known results on transversals. Let $\mathcal{C}$ be a family of finite sets. A \emph{blocking set} of $\mathcal{C}$ is a set that intersects (blocks) every member of $\mathcal{C}$. A blocking set of $\mathcal{C}$ is \emph{minimal} if none of its proper subsets is a blocking set of $\mathcal{C}$.  Minimal blocking sets are also called \emph{transversals}. The family of all minimal blocking sets of $\mathcal{C}$ is denoted $\mathsf{Tr}(\mathcal{C})$ (transversal of $\mathcal{C}$). It can be shown that if $\mathcal{C}_0$ is the family of the inclusion-minimal elements of a family $\mathcal{C}$ of finite sets, then $\mathsf{Tr}(\mathcal{C})=\mathsf{Tr}(\mathcal{C}_0)$ and $\mathsf{Tr}(\mathsf{Tr}(\mathcal{C}))=\mathsf{Tr}(\mathsf{Tr}(\mathcal{C}_0))=\mathcal{C}_0$. (See for instance the book \cite{Ju-2011} for general results on transversals.)
 
\begin{prop}
\label{transversal}
If $\mathcal{F}_i(\mathcal{H})$ and $\mathcal{I}_i(\mathcal{H})$ are the families of inclusion-minimal $\mathcal{R}_i$-forcing sets and $\mathcal{R}_i$-immune sets, respectively, of a hypergraph $\mathcal{H}$, then $\mathsf{Tr}(\mathcal{F}_i(\mathcal{H}))=\mathcal{I}_i(\mathcal{H})$ and $\mathsf{Tr}(\mathcal{I}_i(\mathcal{H}))=\mathcal{F}_i(\mathcal{H})$.
\end{prop}

\begin{pf}
Observe that if $X\subseteq \Omega$ is an $\mathcal{R}_i$-forcing set and if $Y\subseteq \Omega $ is a subset such that $X\cap Y=\emptyset$, then $\Omega\setminus Y$ is an $\mathcal{R}_i$-forcing set, because $X\subseteq \Omega\setminus Y$ and hence ${\cal R}_i^\ast(\Omega\setminus Y)=\Omega$. Since a subset $I\subseteq\Omega$ is by definition an $\mathcal{R}_i$-immune set if and only if $I$ is non-empty and ${\cal R}_i^*(\Omega\setminus I)=\Omega \setminus I$, we deduce that if $X$ is a ${\cal R}_i$-forcing set, then $X\cap I\ne\emptyset$ for all ${\cal R}_i$-immune set $I$. 

Notice also that if $X\subseteq \Omega$ is not an $\mathcal{R}_i$-forcing set, then the set $I=\Omega \setminus \mathcal{R}_i^\ast(X)$ is non-empty and moreover $\mathcal{R}_i^\ast(\Omega\setminus I)=\mathcal{R}_i^\ast(X)=\Omega\setminus I$. We deduce that $I$ is an $\mathcal{R}_i$-immune set and $I\cap\mathcal{R}_i^\ast(X)=\emptyset$. Since $X\subseteq\mathcal{R}_i^\ast(X)$ we conclude that if $X$ is not a  $\mathcal{R}_i$-forcing set, then there exists an $\mathcal{R}_i$-immune set $I$ such that $X\cap I=\emptyset$.

Until now we have proved that if $\Gamma_i=\{I : \, I$ is an $\mathcal{R}_i$-immune set$\}$, then $\mathsf{Tr}(\Gamma_i)=\mathcal{F}_i(\mathcal{H})$ and so  by minimality we have $\mathsf{Tr}(\mathcal{I}_i(\mathcal{H}))=\mathsf{Tr}(\Gamma_i)$ $=\mathcal{F}_i(\mathcal{H})$.  Finally, again by transversality,  $\mathsf{Tr}(\mathcal{F}_i(\mathcal{H}))=\mathsf{Tr}(\mathsf{Tr}(\mathcal{I}_i(\mathcal{H})))=\mathcal{I}_i(\mathcal{H})$.
\end{pf}

\begin{example}
Let $\Omega=\{1,2,3,4\}$ and consider the hypergraph $\mathcal{H}$ on $\Omega$ with hyperedges $\{1,2,3\}$, $\{1,2,4\}$ and $\{1,3,4\}$. We can easily check that the family of minimal forcing sets of ${\cal H}$ are $\mathcal{F}_1(\mathcal{H})=\{\{1,2\},\{1,3\},\{1,4\},\{2,3\},\{2,4\},\{3,4\}\}$ and $\mathcal{F}_2(\mathcal{H})=\{\{2,3\},\{2,4\},\{3,4\}\}$. From the above proposition $\mathcal{I}_i(\mathcal{H})=\mathsf{Tr}(\mathcal{F}_i(\mathcal{H}))$. Hence, by considering the transversal of $\mathcal{F}_i(\mathcal{H})$, $i=1,2$,  we obtain the families of inclusion-minimal immune sets $\mathcal{I}_1(\mathcal{H})=\{\{1,2,3\},\{1,2,4\},$ $\{1,3,4\},\{2,3,4\}\}$ and  $\mathcal{I}_2(\mathcal{H})=\{\{2,3\},\{2,4\},\{3,4\}\}$. 
\end{example}

%%%%%%%%%%%%%%%%%%%%%%%%%%%%%%%%%%%%%%%%%%%%%%%%%%%%%%
\begin{table}[!ht]
\caption{Up to isomorphism, hypergraphs on $\Omega\subseteq \{1,2,3,4\}$ and families of its minimal forcing and minimal immune sets.}
\label{table 1}
{\small 
\begin{center}
\begin{tabular}{|l|c|c|c|c|}
\hline
& & & & \\ 

$\hphantom{aaaaa}{\cal H}_{i,j}$ & 
${\cal F}_1({\cal H}_{i,j})$ & 
${\cal F}_2({\cal H}_{i,j})$ & 
${\cal I}_1({\cal H}_{i,j})$ & 
${\cal I}_2({\cal H}_{i,j})$  \\ 

& & & & \\ 
\hline
& & & & \\ 

${\cal H}_{1,1}=\{\{1\}\}$ & 
${\cal H}_{1,1}$ & 
${\cal H}_{1,1}$ & 
${\cal H}_{1,1}$ & 
${\cal H}_{1,1}$  \\ 

& & & & \\ 
\hline
& & & & \\ 

${\cal H}_{2,1}=\{\{1,2\}\}$ & 
${\cal H}_{2,2}$ & 
${\cal H}_{2,2}$ & 
${\cal H}_{2,1}$ & 
${\cal H}_{2,1}$ \\ 

${\cal H}_{2,2}=\{\{1\},\{2\}\}$ & 
${\cal H}_{2,1}$ & 
${\cal H}_{2,1}$ & 
${\cal H}_{2,2}$ & 
${\cal H}_{2,2}$ \\ 

& & & & \\ 
\hline
& & & & \\ 

${\cal H}_{3,1}=\{\{1,2,3\}\}$ & 
${\cal H}_{3,5}$ & 
${\cal H}_{3,5}$ & 
${\cal H}_{3,1}$ & 
${\cal H}_{3,1}$ \\ 

${\cal H}_{3,2}=\{\{1,2\},\{1,3\},\{2,3\}\}$ & 
${\cal H}_{3,2}$ & 
${\cal H}_{3,2}$ & 
${\cal H}_{3,2}$ & 
${\cal H}_{3,2}$ \\ 

${\cal H}_{3,3}=\{\{1,2\},\{1,3\}\}$ & 
${\cal H}_{2,2}$ & 
${\cal H}_{2,2}$ & 
${\cal H}_{2,1}$ & 
${\cal H}_{2,1}$ \\ 

${\cal H}_{3,4}=\{\{1,2\},\{3\}\} $ & 
${\cal H}_{3,3}$ & 
${\cal H}_{3,3}$ & 
${\cal H}_{3,4}$ & 
${\cal H}_{3,4}$ \\
 
${\cal H}_{3,5}=\{\{1\},\{2\},\{3\}\}$ & 
${\cal H}_{3,1}$ & 
${\cal H}_{3,1}$ & 
${\cal H}_{3,5}$ & 
${\cal H}_{3,5}$ \\

& & & & \\ 
\hline
& & & & \\ 

${\cal H}_{4,1}=\{\{1,2,3,4\}\}$ & 
${\cal H}_{4,20}$ & 
${\cal H}_{4,20}$ &  
${\cal H}_{4,1}$ & 
${\cal H}_{4,1}$ \\ 

${\cal H}_{4,2}=\{\{1,2,3\},\{1,2,4\},\{1,3,4\},\{2,3,4\}\}$
& ${\cal H}_{4,10}$ & 
${\cal H}_{4,2}$ & 
${\cal H}_{4,2}$ & 
${\cal H}_{4,10}$ \\

${\cal H}_{4,3}=\{\{1,2,3\},\{1,2,4\},\{1,3,4\}\}$ & 
${\cal H}_{4,10}$ & 
${\cal H}_{3,2}$ & 
${\cal H}_{4,2}$ & 
${\cal H}_{3,2}$ \\ 

${\cal H}_{4,4}=\{\{1,2,3\},\{1,2,4\}\}$ & 
${\cal H}_{2,2}$ & 
${\cal H}_{2,2}$ & 
${\cal H}_{2,1}$ & 
${\cal H}_{2,1}$ \\ 

${\cal H}_{4,5}=\{\{1,2,3\},\{1,2,4\},\{3,4\}\}$ & 
${\cal H}_{4,11}$ & 
${\cal H}_{4,11}$ & 
${\cal H}_{4,5}$ & 
${\cal H}_{4,5}$ \\ 

${\cal H}_{4,6}=\{\{1,2,3\},\{1,4\}\}$ & 
${\cal H}_{3,5}$ & 
${\cal H}_{3,5}$ & 
${\cal H}_{3,1}$ & 
${\cal H}_{3,1}$ \\ 

${\cal H}_{4,7}=\{\{1,2,3\},\{1,4\},\{2,4\}\}$ & 
${\cal H}_{4,16}$ & 
${\cal H}_{4,16}$ & 
${\cal H}_{4,3}$ & 
${\cal H}_{4,3}$ \\ 

${\cal H}_{4,8}=\{\{1,2,3\},\{1,4\},\{2,4\},\{3,4\}\}$ & 
${\cal H}_{4,10}$ & 
${\cal H}_{4,10}$ & 
${\cal H}_{4,2}$ & 
${\cal H}_{4,2}$ \\ 

${\cal H}_{4,9}=\{\{1,2,3\},\{4\}\}$ & 
${\cal H}_{4,14}$ & 
${\cal H}_{4,14}$ & 
${\cal H}_{4,9}$ & 
${\cal H}_{4,9}$ \\ 

${\cal H}_{4,10}=\{\{1,2\},\{1,3\},\{1,4\},\{2,3\},\{2,4\},\{3,4\}\}$ & 
${\cal H}_{4,2}$ & \
${\cal H}_{4,2}$ & 
${\cal H}_{4,10}$ & 
${\cal H}_{4,10}$ \\ 

${\cal H}_{4,11}=\{\{1,2\},\{1,3\},\{1,4\},\{2,3\},\{2,4\}\}$ & 
${\cal H}_{4,13}$ & 
${\cal H}_{4,13}$ & 
${\cal H}_{4,17}$ & 
${\cal H}_{4,17}$ \\ 

${\cal H}_{4,12}=\{\{1,2\},\{1,3\},\{1,4\},\{2,3\}\}$ & 
${\cal H}_{4,11}$ & 
${\cal H}_{4,11}$ & 
${\cal H}_{4,5}$ & 
${\cal H}_{4,5}$ \\ 

${\cal H}_{4,13}=\{\{1,2\},\{1,3\},\{2,4\},\{3,4\}\}$ & 
${\cal H}_{4,13}$ & 
${\cal H}_{4,13}$ & 
${\cal H}_{4,17}$ & 
${\cal H}_{4,17}$ \\ 

${\cal H}_{4,14}=\{\{1,2\},\{1,3\},\{1,4\}\}$ & 
${\cal H}_{3,2}$ & 
${\cal H}_{3,2}$ & 
${\cal H}_{3,2}$ & 
${\cal H}_{3,2}$ \\ 

${\cal H}_{4,15}=\{\{1,2\},\{1,3\},\{2,4\}\}$ & 
${\cal H}_{4,19}$ & 
${\cal H}_{4,19}$ & 
${\cal H}_{4,4}$ & 
${\cal H}_{4,4}$ \\ 

${\cal H}_{4,16}=\{\{1,2\},\{1,3\},\{2,3\},\{4\}\}$ & 
${\cal H}_{4,3}$ & 
${\cal H}_{4,3}$ & 
${\cal H}_{4,16}$ & 
${\cal H}_{4,16}$ \\ 

${\cal H}_{4,17}=\{\{1,2\},\{3,4\}\}$ & 
${\cal H}_{4,13}$ & 
${\cal H}_{4,13}$ & 
${\cal H}_{4,17}$ & 
${\cal H}_{4,17}$ \\ 

${\cal H}_{4,18}=\{\{1,2\},\{1,3\},\{4\}\}$ &
${\cal H}_{3,3}$ & 
${\cal H}_{3,3}$ & 
${\cal H}_{3,4}$ & 
${\cal H}_{3,4}$ \\ 

${\cal H}_{4,19}=\{\{1,2\},\{3\},\{4\}\}$ & 
${\cal H}_{4,4}$ & 
${\cal H}_{4,4}$ & 
${\cal H}_{4,19}$ & 
${\cal H}_{4,19}$ \\ 

${\cal H}_{4,20}=\{\{1\},\{2\},\{3\},\{4\}\}$ & 
${\cal H}_{4,1}$ & 
${\cal H}_{4,1}$ & 
${\cal H}_{4,20}$ & 
${\cal H}_{4,20}$ \\ 
& & & & \\ 
\hline
\end{tabular}
\end{center}
}
\end{table}
%%%%%%%%%%%%%%%%%%%%%%%%%%%%%%%%%%%%%%%%%%%%%%%%%%%

\begin{example}
If $\Omega\subseteq \{1,\ldots, i\}$, $1\leqslant i\leqslant 4$, then, up to isomorphism, there exist 28 hypergraphs ${\cal H}_{i,j}$ which are enumerated in Table~\ref{table 1}. For simplicity, each hypergraph $\mathcal{H}_{i,j}$ is identified with its hyperedge clutter, and for each $\mathcal{H}_{i,j}$ the clutters $\mathcal{F}_1(\mathcal{H}_{i,j})$, $\mathcal{F}_2(\mathcal{H}_{i,j})$, $\mathcal{I}_1(\mathcal{H}_{i,j})$ and $\mathcal{I}_2(\mathcal{H}_{i,j})$ have been explicitly calculated and displayed (up to isomorphism). For instance, if ${\cal H}_{3,3}=\{\{1,2\},\{2,3\}\}$, then ${\cal F}_1({\cal H}_{3,3})={\cal F}_2({\cal H}_{3,3})=\{\{1\},\{3\}\}\cong{\cal H}_{2,2}$ and ${\cal I}_1({\cal H}_{3,3})={\cal I}_2({\cal H}_{3,3})=\{\{1,3\}\}\cong{\cal H}_{2,1}$; and if ${\cal H}_{3,4}=\{\{1,2\},\{3\}\}$, then ${\cal F}_1({\cal H}_{3,4})={\cal F}_2({\cal H}_{3,4})=\{\{1,3\},\{2,3\}\}\cong{\cal H}_{3,3}$ and ${\cal I}_1({\cal H}_{3,4})={\cal I}_2({\cal H}_{3,4})=\{\{1,2\},\{3\}\}\cong{\cal H}_{3,4}$ (in this last case we have actually an equality).
\end{example}

From the analysis of Table~\ref{table 1} we deduce the following fact: if ${\cal H}$ is a hypergraph on a finite set $\Omega$ with cardinality $|\Omega|\leqslant 4$, then, up to isomorphism, ${\cal F}_1({\cal H})= {\cal F}_2({\cal H})$ if and only if ${\cal I}_1({\cal H})={\cal I}_2({\cal H})$ if and only if  ${\cal H}$ is not isomorphic to ${\cal H}_{4,2}$, neither ${\cal H}$ is not isomorphic to ${\cal H}_{4,3}$.

The following proposition is the only general result that we provide concerning the relation between forcing sets and immune sets with respect to the rules $\mathcal{R}_1$ and $\mathcal{R}_2$.

\begin{prop}
\label{R1R2}
Let $\mathcal{H}$ be a hypergraph on a finite set $\Omega$. Then the following statements hold:
\begin{enumerate}
\item If $F\subseteq\Omega$ is an $\mathcal{R}_2$-forcing set, then $F$ is also an $\mathcal{R}_1$-forcing set. Moreover, if $F\in\mathcal{F}_2(\mathcal{H})$ is a minimal $\mathcal{R}_2$-forcing set, then there exists a minimal $\mathcal{R}_1$-forcing set $F'\in\mathcal{F}_1(\mathcal{H})$ such that  $F'\subseteq F$.
\item If $I\subseteq\Omega$ is an $\mathcal{R}_1$-immune set, then $I$ is also an $\mathcal{R}_2$-immune set. Moreover, if $I\in\mathcal{I}_1(\mathcal{H})$ is a minimal $\mathcal{R}_1$-immune set, then there exists a minimal $\mathcal{R}_2$-immune set $I'\in\mathcal{I}_2(\mathcal{H})$ such that  $I'\subseteq I$.
\end{enumerate}
\end{prop}

\begin{pf}
Firstly, let us demonstrate the following claim: if $B$ is the set of vertices of $\mathcal{H}$ which are initially colored black, then $\mathcal{R}^\ast_2(B) \subseteq \mathcal{R}^\ast_1(B)$. Let us prove our claim.  Clearly we have $B \subseteq \mathcal{R}^\ast_1(B)$ and hence $\mathcal{R}^\ast_2(B) \subseteq \mathcal{R}^\ast_2(\mathcal{R}^\ast_1(B))$, because the rule $\mathcal{R}_2$ is increasing. Let $B'=\mathcal{R}^\ast_1(B)$. We are done if we demonstrate that $\mathcal{R}^\ast_2(B')=B'$. In agreement with the definitions of $\mathcal{R}_1$ and $\mathcal{R}_2$ we have that $\mathcal{R}^\ast_2(B')\subseteq \mathcal{R}^\ast_1(B')$. So we have that $B'\subseteq \mathcal{R}^\ast_2(B')\subseteq \mathcal{R}^\ast_1(B')=B'$. Therefore we conclude that  $\mathcal{R}^\ast_2(B')=B'$.

Now let us prove the proposition. By our claim we have that $\mathcal{R}^\ast_2(F)\subseteq \mathcal{R}^\ast_1(F)$ and $\mathcal{R}^\ast_2(V(\mathcal{H})\setminus I)\subseteq \mathcal{R}^\ast_1(V(\mathcal{H})\setminus I)$. Therefore any $\mathcal{R}_2$-forcing set $F$ is also an $\mathcal{R}_1$-forcing set and any $\mathcal{R}_1$-immune set $I$ is also an $\mathcal{R}_2$-immune set. To finish the proof of the proposition observe  that the ``moreover'' condition of both statements is a basic and general result. Namely, if $ \Gamma_1$ and $\Gamma_2$ are two families of subsets of a set $X$ such that $\Gamma_1\subseteq \Gamma_2$, then for all inclusion-minimal element $Y$ of $\Gamma_1$ there exists an inclusion-minimal element $Z$ of $\Gamma_2$ such that $Z \subseteq  Y$.
\end{pf}

%%%%%%%%%%%%%%%%%%%%%%%%%%%%%%%%%%%%%%%%%%%%%%%%%%%%%%%%
\section{Neighbourhood characterization of immune sets.}
\label{Neighbourhood characterization}	
For zero forcing, a nice characterization of the immune sets  can be given in terms of neighbourhoods, without making reference to the iterative zero forcing process. Let us consider first the case in which the hypergraph is a graph $G$. Given $X\subseteq V(G)$, let $V(G)\setminus X$ be the initial set of black vertices. Observe that a white vertex of $X$ can be forced (by a black vertex $v\in V(G)\setminus X$) to become black if and only if $|N(v)\cap X|=1$. So we have the following fact.

\begin{rmk}
\label{graph-immune}
Let $G$ be a graph. A non-empty subset  $X\subseteq V(G)$ is an immune set of $G$ if and only if $|N(v)\cap X|\ne 1$ for all $v\in V(G)\setminus X$.
\end{rmk}

Let $\mathcal{H}$  be a hypergraph on a finite set $\Omega$. Given a subset $B\subseteq \Omega$, let us define the open neighbourhood of $B$, $\mathcal{N}(B)$, as $\mathcal{N}(B)=\{B'\subseteq \Omega\, : \, B\cap B'=\emptyset$ and $B\cup B'\in E({\cal H})\}$. Notice that if $\mathcal{H}$ is a graph $G$ and $v$ is a vertex, then $\mathcal{N}(\{v\})=\{\{w\} \, : \, w\in N(v)\}$. Hence if $G$ is seen as a hypergraph $\mathcal{H}$, Remark~\ref{graph-immune} can be alternatively formulated as follows:

\begin{prop}
\label{graph-immune-bis}
Let $G$ be a graph. A non-empty subset  $X\subseteq V(G)$ is an immune set of $G$ if and only if $|\{B'\in {\cal N}(\{v\}) \, : \, B'\cap X\neq \emptyset \}|\ne 1$ for all $\{v\}\subseteq V(G)\setminus X$.
\end{prop}

If in a hypergraph $\mathcal{H}$ we consider the zero forcing rule $\mathcal{R}_2$, then Proposition~\ref{graph-immune-bis} can be straightforwardly generalized.
\begin{prop}
\label{hypergraph-immune2}
Let $\mathcal{H}$  be a hypergraph on a finite set $\Omega$. A non-empty subset  $X\subseteq \Omega$ is an $\mathcal{R}_2$-immune set of ${\cal H}$  if and only if $|\{B'\in \mathcal{N}(B) \, : \, B'\cap X\neq \emptyset \}|\ne 1$ for all $B\subseteq \Omega\setminus X$.
\end{prop}

\begin{pf}
Given $\emptyset\ne X\subseteq \Omega$, let $\Omega\setminus X$ be the initial set of black vertices and  let $A$ be a hyperedge such that $A\setminus X\ne\emptyset$ and $A\cap X\ne\emptyset$. According to the zero forcing rule $\mathcal{R}_2$ the black subset $A\setminus X$ of the hyperedge $A$ will force the white subset $A\cap X$ to become black if and only if there is not another hyperedge with white vertices and containing $A\setminus X$. Hence a black subset $B\subseteq A\setminus X\subseteq \Omega\setminus X$  of an hyperedge $A$ will force $A\cap X$ to become black if and only if $|\{B'\in \mathcal{N}(B) \, : \, B'\cap X\neq \emptyset \}|= 1$. 
\end{pf}

Observe that if we use the standard definition in hypergraphs of the neighbourhood of a subset of vertices $B$, namely $N(B)=\{A\in E(\mathcal{H}): B\subseteq A\}$, then the necessary and sufficient condition in the above proposition can also be expressed as: $|\{ A'\in N(A\setminus X) \, : \,  A'\cap X\neq \emptyset\}|\geqslant 2$ for all hyperedge $A$ such that $A\setminus X\ne\emptyset$ and $A\cap X\ne\emptyset$. This observation motivates the following notation, that will allow us to characterize also the $\mathcal{R}_1$-immune sets of $\mathcal{H}$. Let ${\cal H}$ be a hypergraph on a finite set $\Omega$. Let $X\subseteq \Omega$ be a subset of vertices, and let $A\in E({\cal H})$ be a hyperedge. Let us define the following sets:
\begin{align*}
\Sigma_1(X,A) &
%=\{ A'\in E({\cal H}) \, : \, A\setminus X\subseteq A' \textrm{ and } (A'\cap X) \setminus A \neq  \emptyset \}
=\{ A'\in N(A\setminus X) \, : \,  (A'\cap X) \setminus A \neq  \emptyset \};\\
\Sigma_2(X,A) &
%=\{ A'\in E({\cal H}) \, : \, A\setminus X\subseteq A' \textrm{ and } A'\cap X\neq \emptyset\}.
=\{ A'\in N(A\setminus X) \, : \,  A'\cap X\neq \emptyset\}.
\end{align*}
Notice that if $A$ is a hyperedge  of ${\cal H}$ such that $A\setminus X \neq \emptyset$ and $A\cap X\neq \emptyset$, then
we have  $\Sigma_1(X,A)\cup \{A\} \subseteq\Sigma_2(X,A)$ but,  in general, the above inclusion is not an equality.

Using these sets we have the following characterizations of immune sets in graphs and hypergraphs.

\begin{prop}
\label{graph-immune-bis-2}
Let $G$ be a graph, and let  $X\subseteq V(G)$ be a non-empty subset. Then the following conditions are equivalent:
\begin{enumerate}
\item $X$ is an immune set of $G$ with respect to the zero forcing rule.
\item $|\Sigma_1(X,e)|\geqslant 1$ for all edges $e\in E(G)$ such that $e\setminus X \neq \emptyset$ and $e\cap X\neq \emptyset$.
\item  $|\Sigma_2(X,e)|\geqslant 2$ for all edges $e\in E(G)$ such that $e\setminus X \neq \emptyset$ and $e\cap X\neq \emptyset$.
\end{enumerate}
\end{prop}

\begin{pf}
The result follows from Remark~\ref{graph-immune} and the fact that if $G$ is a graph and $e=\{v,w\}\in E(G)$ is an edge such that $e\setminus X=\{v\}$ and $e\cap X=\{w\}$,  then $\Sigma_2(X,e)=\Sigma_1(X,e)\cup \{e\}$ and $|\Sigma_2(X,e)|\geqslant 2$ if and only if $|\Sigma_1(X,e)|\geqslant 1$, if and only if $|N(v)\cap X|\geqslant 2$.
\end{pf}

\begin{prop}
\label{hyper-immune}
Let ${\cal H}$ be a hypergraph  on a finite set $\Omega$ and let $X\subseteq \Omega $ be non-empty subset. Then the following statements hold:
\begin{enumerate}
\item $X$ is an $\mathcal{R}_1$-immune set if and only if $|\Sigma_1(X,A)|\geqslant 1$ for all hyperedges $A$ such that $A\setminus X\ne\emptyset$ and $A\cap X\ne\emptyset$. 
\item $X$ is an $\mathcal{R}_2$-immune set if and only if $|\Sigma_2(X,A)|\geqslant 2$ for all hyperedges $A$ such that $A\setminus X\ne\emptyset$ and $A\cap X\ne\emptyset$. 
\end{enumerate}
\end{prop}

\begin{pf}
Given $\emptyset\ne X\subseteq V(\mathcal{H})$, let $V(\mathcal{H})\setminus X$ be the initial set of black vertices and  let $A$ be a hyperedge such that $A\setminus X\ne\emptyset$ and $A\cap X\ne\emptyset$. 

According to the zero forcing rule $\mathcal{R}_1$ the black subset $A\setminus X$ of the hyperedge $A$ will force the white subset $A\cap X$ to become black if and only if there are no white vertices outside of $A$ adjacent to $A\setminus X$. Therefore $A\setminus X$ will force $A\cap X$ to become black if and only if the set of white vertices adjacent to $A\setminus X$ is contained in $A\cap X$; if and only if $\Sigma_1(X,A)=\emptyset$. This proves the part of the proposition concerning $\mathcal{R}_1$.

Analogously, according to the zero forcing rule $\mathcal{R}_2$ the black subset $A\setminus X$ of the hyperedge $A$ will force the white subset $A\cap X$ to become black if and only if there is not another hyperedge $A'$ with white vertices and containing $A\setminus X$; if and only if $\Sigma_2(X,A)=\{A\}$. This finishes the proof of the proposition.
\end{pf}

\begin{rmk}
Notice that if $|\Sigma_1(X,A)|\geqslant 1$, then $|\Sigma_2(X,A)|\geqslant 2$ and so if $X$ is an $\mathcal{R}_1$-immune set, then it is also an $\mathcal{R}_2$-immune set, as stated in Proposition~\ref{R1R2}.
\end{rmk}

%%%%%%%%%%%%%%%%%%%%%%%%%%%%%%%%%%%%%%%%%%%%%%%%%%%%%%%%%%%%%
\section{Characterization of uniform forcing and immune clutters.}
\label{clutters} 

As discussed in Section~\ref{ZeroForcing}, given a hypergraph $\mathcal{H}$, the families $\mathcal{F}_i(\mathcal{H})$ and $\mathcal{I}_i(\mathcal{H})$ of its (inclusion-)minimal forcing and immune sets with respect to the zero forcing rule $\mathcal{R}_i$ constitute a pair of clutters,  each one being the transversal of the other. (If the hypergraph is a graph $G$, we can emphasize this fact by denoting its families of minimal forcing and minimal immune sets by $\mathcal{F}(G)$ and $\mathcal{I}(G)$, respectively.) Two interesting open problems arise at this point: given a clutter $\Delta$ on a finite set $\Omega$, is there a graph $G$ or, more generally, a hypergraph $\mathcal{H}$, with vertex set $\Omega$, such that its family of minimal forcing or immune sets with respect to the zero forcing rule $\mathcal{R}_i$ is precisely $\Delta$; and if the answer is in the affirmative, can we find $G$ or $\mathcal{H}$? These questions motivate the following definition. 

\begin{defi} 
\label{defi-reali}
Let $\Delta$ be a clutter on a finite set $\Omega$. Then
\begin{enumerate}
\item[(a)] $\Delta$ is  a graph-forcing clutter (respectively, an $\mathcal{R}_i$-forcing clutter) if there exists a graph $G$ (respectively, a hypergraph $\mathcal{H}$) such that $V(G)=\Omega$ and ${\cal F}(G)=\Delta$ (respectively, $V(\mathcal{H})=\Omega$ and ${\cal F}_i(\mathcal{H})=\Delta$). In this case, we say that $G$ is a graph-forcing realization (respectively, $\mathcal{H}$ is an $\mathcal{R}_i$-forcing realization) of $\Delta$.
\item[(b)]  $\Delta$ is a graph-immune clutter (respectively, an $\mathcal{R}_i$-immune clutter) if there exists a graph $G$ (respectively, a hypergraph $\mathcal{H}$) such that $V(G)=\Omega$ and ${\cal I}(G)=\Delta$ (respectively, $V(\mathcal{H})=\Omega$ and ${\cal I}_i(\mathcal{H})=\Delta$). In this case, we say that $G$ is a graph-immune realization (respectively, $\mathcal{H}$ is an $\mathcal{R}_i$-immune realization) of $\Delta$.		
\end{enumerate}
\end{defi}

To illustrate Definition~\ref{defi-reali} let us consider the following example:

\begin{example}
From Table~\ref{table 1} we also deduce which ones of the $28$ non-isomorphic clutters $\Delta$ defined on a set $\Omega$ with cardinality $|\Omega|\leqslant 4$ are $\mathcal{R}_i$-forcing clutters and which ones are $\mathcal{R}_i$-immune clutters. This analysis is detailed in Table~\ref{table 2}. We underline the following facts:
\begin{enumerate}
\item
$\Delta$ is an $R_1$-forcing clutter (respectively, $\Delta$ is an $R_1$-immune  clutter) if and only if $\Delta$ is an $R_2$-forcing clutter (respectively, $\Delta$ is an $R_2$-immune clutter).

\item There are cases in which $\Delta$ is both an $R_i$-forcing  clutter and an $R_i$-immune clutter; for example, $\mathcal{H}_{4,16}=\mathcal{F}_1(\mathcal{H}_{4,7})=\mathcal{I}_1(\mathcal{H}_{4,16})$. There are samples in which $\Delta$ is an $R_i$-forcing  clutter but not an $R_i$-immune clutter; for instance, $\mathcal{H}_{4,13}=\mathcal{F}_2(\mathcal{H}_{4,11})$ but there is no realization of $\mathcal{H}_{4,13}$ as an $\mathcal{R}_2$-immune clutter.  There are also cases in which $\Delta$ is an $R_i$-immune clutter but not an $R_i$-forcing clutter. Finally, observe that $\Delta$ is not an $R_i$-forcing clutter nor an $R_i$-immune clutter if and only if $\Delta\cong {\cal H}_{4,6}, {\cal H}_{4,7}$,  ${\cal H}_{4,8}, {\cal H}_{4,12},{\cal H}_{4,15},{\cal H}_{4,18}$.

\item The realizations of a given $\Delta$ may depend  on $\mathcal{R}_i$  and, moreover, the number of non-isomorphic realizations  may be different depending on $\Delta$. 
\end{enumerate}
 
%%%%%%%%%%%%%%%%%%%%%%%%%%%%%%%%%%%%%%%%%%%%%%%%%%%%%%%%%%%%%%%%%%%%
\begin{table}[!ht]
\caption{Up to isomorphism, clutters $\Delta$ on $\Omega\subseteq \{1,2,3,4\}$ and realizations of $\Delta$ as $\mathcal{R}_i$-forcing and $\mathcal{R}_i$-immune clutters.}
\label{table 2}
{\small 
\begin{center}
\begin{tabular}{|c|c|c|c|c|}
\hline
& & & & \\ 

& 
${\cal H}_{k,\ell}$ such that & 
${\cal H}_{k,\ell}$ such that & 
${\cal H}_{k,\ell}$ such that & 
${\cal H}_{k,\ell}$ such that  \\

$\Delta ={\cal H}_{i,j} $ & 
${\cal F}_1({\cal H}_{k,\ell})=\Delta $ & 
${\cal F}_2({\cal H}_{k,\ell})=\Delta $ & 
${\cal I}_1({\cal H}_{k,\ell})=\Delta $ & 
${\cal I}_2   ({\cal H}_{k,\ell})=\Delta $ \\ 

& & & & \\ 
\hline
& & & & \\ 

${\cal H}_{1,1}$ & 
${\cal H}_{1,1}$ & 
${\cal H}_{1,1}$ & 
${\cal H}_{1,1}$ & 
${\cal H}_{1,1}$  \\ 

& & & & \\ 
\hline
& & & & \\ 

${\cal H}_{2,1}$ & 
${\cal H}_{2,2}$ & 
${\cal H}_{2,2}$ & 
${\cal H}_{2,1}, {\cal H}_{3,3}, {\cal H}_{4,4}$ & 
${\cal H}_{2,1}, {\cal H}_{3,3}, {\cal H}_{4,4}$ \\ 

${\cal H}_{2,2}$ & 
${\cal H}_{2,1}, {\cal H}_{3,3}, {\cal H}_{4,4}$ & 
${\cal H}_{2,1}, {\cal H}_{3,3}, {\cal H}_{4,4}$ &
${\cal H}_{2,2}$ & 
${\cal H}_{2,2}$ \\ 

& & & & \\ 
\hline
& & & & \\ 

${\cal H}_{3,1}$ & 
${\cal H}_{3,5}$ & 
${\cal H}_{3,5}$ & 
${\cal H}_{3,1},{\cal H}_{4,6}$ & 
${\cal H}_{3,1},{\cal H}_{4,6}$ \\ 

${\cal H}_{3,2}$ & 
${\cal H}_{3,2}, {\cal H}_{4,14}$ & 
${\cal H}_{3,2}, {\cal H}_{4,3}, {\cal H}_{4,14}$ & 
${\cal H}_{3,2}, {\cal H}_{4,14}$ & 
${\cal H}_{3,2}, {\cal H}_{4,3}, {\cal H}_{4,14}$ \\ 

${\cal H}_{3,3}$ & 
${\cal H}_{3,4}, {\cal H}_{4,18}$ & 
${\cal H}_{3,4}, {\cal H}_{4,18}$ & 
$\hphantom{aaa}$ \hrulefill $\hphantom{aaa}$ & 
$\hphantom{aaa}$ \hrulefill $\hphantom{aaa}$ \\ 

${\cal H}_{3,4}$ & 
$\hphantom{aaa}$ \hrulefill $\hphantom{aaa}$ & 
$\hphantom{aaa}$ \hrulefill $\hphantom{aaa}$ & 
${\cal H}_{3,4}, {\cal H}_{4,18}$ & 
${\cal H}_{3,4}, {\cal H}_{4,18}$ \\
 
${\cal H}_{3,5}$ & 
${\cal H}_{3,1}, {\cal H}_{4,6}$ & 
${\cal H}_{3,1}, {\cal H}_{4,6}$ & 
${\cal H}_{3,5}$ & 
${\cal H}_{3,5}$ \\

& & & & \\ 
\hline
& & & & \\ 

${\cal H}_{4,1}$ & 
${\cal H}_{4,20}$ & 
${\cal H}_{4,20}$ &  
${\cal H}_{4,1}$ & 
${\cal H}_{4,1}$ \\ 

${\cal H}_{4,2}$ & 
${\cal H}_{4,10}$ & 
${\cal H}_{4,2}, {\cal H}_{4,10}$ & 
${\cal H}_{4,2},{\cal H}_{4,3}, {\cal H}_{4,8}$ & 
${\cal H}_{4,8}$ \\

${\cal H}_{4,3}$ & 
${\cal H}_{4,16}$ & 
${\cal H}_{4,16}$ & 
${\cal H}_{4,7}$ & 
${\cal H}_{4,7}$ \\ 

${\cal H}_{4,4}$ & 
${\cal H}_{4,19}$ & 
${\cal H}_{4,19}$ & 
${\cal H}_{4,15}$ & 
${\cal H}_{4,15}$ \\ 

${\cal H}_{4,5}$ & 
$\hphantom{aaa}$ \hrulefill $\hphantom{aaa}$ & 
$\hphantom{aaa}$ \hrulefill $\hphantom{aaa}$ & 
${\cal H}_{4,5}, {\cal H}_{4,12}$ & 
${\cal H}_{4,5}, {\cal H}_{4,12}$ \\ 

${\cal H}_{4,6}$ & 
$\hphantom{aaa}$ \hrulefill $\hphantom{aaa}$ & 
$\hphantom{aaa}$ \hrulefill $\hphantom{aaa}$ & 
$\hphantom{aaa}$ \hrulefill $\hphantom{aaa}$ & 
$\hphantom{aaa}$ \hrulefill $\hphantom{aaa}$ \\ 

${\cal H}_{4,7}$ & 
$\hphantom{aaa}$ \hrulefill $\hphantom{aaa}$ & 
$\hphantom{aaa}$ \hrulefill $\hphantom{aaa}$ & 
$\hphantom{aaa}$ \hrulefill $\hphantom{aaa}$ & 
$\hphantom{aaa}$ \hrulefill $\hphantom{aaa}$ \\ 

${\cal H}_{4,8}$ & 
$\hphantom{aaa}$ \hrulefill $\hphantom{aaa}$ & 
$\hphantom{aaa}$ \hrulefill $\hphantom{aaa}$ & 
$\hphantom{aaa}$ \hrulefill $\hphantom{aaa}$ & 
$\hphantom{aaa}$ \hrulefill $\hphantom{aaa}$ \\ 

${\cal H}_{4,9}$ & 
$\hphantom{aaa}$ \hrulefill $\hphantom{aaa}$ & 
$\hphantom{aaa}$ \hrulefill $\hphantom{aaa}$ & 
${\cal H}_{4,9}$ & 
${\cal H}_{4,9}$ \\ 

${\cal H}_{4,10}$ & 
${\cal H}_{4,2}, {\cal H}_{4,3}, {\cal H}_{4,8}$ & 
${\cal H}_{4,8}$ & 
${\cal H}_{4,10}$ & 
${\cal H}_{4,2}, {\cal H}_{4,10}$ \\ 

${\cal H}_{4,11}$ & 
${\cal H}_{4,5},{\cal H}_{4,12}$ & 
${\cal H}_{4,5},{\cal H}_{4,12}$ & 
$\hphantom{aaa}$ \hrulefill $\hphantom{aaa}$ & 
$\hphantom{aaa}$ \hrulefill $\hphantom{aaa}$ \\ 

${\cal H}_{4,12}$ &  
$\hphantom{aaa}$ \hrulefill $\hphantom{aaa}$ & 
$\hphantom{aaa}$ \hrulefill $\hphantom{aaa}$ & 
$\hphantom{aaa}$ \hrulefill $\hphantom{aaa}$ & 
$\hphantom{aaa}$ \hrulefill $\hphantom{aaa}$ \\ 

${\cal H}_{4,13}$ & 
${\cal H}_{4,11}, {\cal H}_{4,13}, {\cal H}_{4,17}$ & 
${\cal H}_{4,11}, {\cal H}_{4,13}, {\cal H}_{4,17}$ & 
$\hphantom{aaa}$ \hrulefill $\hphantom{aaa}$ & 
$\hphantom{aaa}$ \hrulefill $\hphantom{aaa}$ \\ 

${\cal H}_{4,14}$ & 
${\cal H}_{4,9}$ & 
${\cal H}_{4,9}$ & 
$\hphantom{aaa}$ \hrulefill $\hphantom{aaa}$ & 
$\hphantom{aaa}$ \hrulefill $\hphantom{aaa}$ \\ 

${\cal H}_{4,15}$ & 
$\hphantom{aaa}$ \hrulefill $\hphantom{aaa}$ & 
$\hphantom{aaa}$ \hrulefill $\hphantom{aaa}$ & 
$\hphantom{aaa}$ \hrulefill $\hphantom{aaa}$ & 
$\hphantom{aaa}$ \hrulefill $\hphantom{aaa}$ \\ 

${\cal H}_{4,16}$ & 
${\cal H}_{4,7}$ & 
${\cal H}_{4,7}$ & 
${\cal H}_{4,16}$ & 
${\cal H}_{4,16}$ \\ 

${\cal H}_{4,17}$ & 
$\hphantom{aaa}$ \hrulefill $\hphantom{aaa}$ & 
$\hphantom{aaa}$ \hrulefill $\hphantom{aaa}$ & 
${\cal H}_{4,11}, {\cal H}_{4,13}, {\cal H}_{4,17}$ & 
${\cal H}_{4,11}, {\cal H}_{4,13}, {\cal H}_{4,17}$ \\ 

${\cal H}_{4,18}$ & 
$\hphantom{aaa}$ \hrulefill $\hphantom{aaa}$ & 
$\hphantom{aaa}$ \hrulefill $\hphantom{aaa}$ & 
$\hphantom{aaa}$ \hrulefill $\hphantom{aaa}$ & 
$\hphantom{aaa}$ \hrulefill $\hphantom{aaa}$ \\ 

${\cal H}_{4,19}$ & 
${\cal H}_{4,15}$ & 
${\cal H}_{4,15}$ & 
${\cal H}_{4,19}$ & 
${\cal H}_{4,19}$ \\ 

${\cal H}_{4,20}$ & 
${\cal H}_{4,1}$ & 
${\cal H}_{4,1}$ & 
${\cal H}_{4,20}$ & 
${\cal H}_{4,20}$ \\ 

& & & & \\ 
\hline
\end{tabular}
\end{center}}
\end{table}

\end{example}
%%%%%%%%%%%%%%%%%%%%%%%%%%%%%%%%%%%%%%%%%%%%%%%%%%%%%

Our goal in this section is to provide a complete characterization of forcing and immune uniform clutters. The $k$-\emph{uniform} clutter ${\cal U}_{k, \Omega}$ on a finite set $\Omega$ is ${\cal U}_{k, \Omega}=\{A \subseteq \Omega\, : \, |A|=k\}$.

The transversal of ${\cal U}_{k, \Omega}$ is also a uniform clutter on $\Omega$; namely, $\mathsf{Tr}({\cal U}_{k, \Omega})={\cal U}_{|\Omega|-k+1, \Omega}$. Thus, by Proposition~\ref{transversal}, the determination of forcing and immune uniform clutters are equivalent problems. 

In the first place let us obtain the families $\mathcal{F}_i(\mathcal{H})$ and $\mathcal{I}_i(\mathcal{H})$ for a \emph{complete hypergraph} on a finite set $\Omega$; that is, a hypergraph $\mathcal{H}_{k,\Omega}$ whose collection of hyperedges is the uniform clutter $\mathcal{U}_{k,\Omega}$. 

Observe that if $\Omega$ is a finite set of size $n$, then  the complete hypergraph $\mathcal{H}_{2,\Omega}$ is just the complete graph $K_n$. If we consider zero forcing, $\mathcal{R}_0$, on $K_n$ it is easily checked that the minimal forcing sets are the sets of cardinality $n-1$, whereas the minimal immune sets are the ones of cardinality two. Therefore we have the following fact.

\begin{prop}
\label{K_n}
Let $K_n$ be the complete graph of $n$ vertices. Then $\mathcal{F}(K_n)=\mathcal{U}_{n-1,\Omega}$ and $\mathcal{I}(K_n)=\mathcal{U}_{2,\Omega}$.
\end{prop}

Next we generalise the above proposition to hypergraphs. Namely, we compute the $\mathcal{R}_i$-forcing sets and the $\mathcal{R}_i$-immune sets of  the complete hypergraph $\mathcal{H}_{k,\Omega}$. Our results are stated in the following two propositions that highlight the difference between the rules $\mathcal{R}_1$ and $\mathcal{R}_2$. 

\begin{prop}
\label{uni-R2}
Let $\mathcal{H}_{k,\Omega}$ be a complete hypergraph on a finite set $\Omega$, $1\leqslant k\leqslant |\Omega|$.  Then the following statements hold:
\begin{enumerate}
\item $\mathcal{F}_2(\mathcal{H}_{1,\Omega})=\mathcal{U}_{|\Omega|,\Omega}$ and $\mathcal{I}_2(\mathcal{H}_{1,\Omega})=\mathcal{U}_{1,\Omega}$.	
\item $\mathcal{F}_2(\mathcal{H}_{k,\Omega})=\mathcal{U}_{|\Omega|-1,\Omega}$ and $\mathcal{I}_2(\mathcal{H}_{k,\Omega})=\mathcal{U}_{2,\Omega}$ for $2\leqslant k\leqslant |\Omega|-1$.
\item $\mathcal{F}_2(\mathcal{H}_{|\Omega|,\Omega})=\mathcal{U}_{1,\Omega}$ and $\mathcal{I}_2(\mathcal{H}_{|\Omega|,\Omega})=\mathcal{U}_{|\Omega|,\Omega}$.
\end{enumerate}
\end{prop}

\begin{pf}
Let $\Omega=\{1,\ldots,n\}$. First observe that the set of hyperedges of  $\mathcal{H}_{1,\Omega}$ is $\{\{1\},\ldots,\{n\}\}$ and so for $i=1,2$ we have $\mathcal{F}_i(\mathcal{H}_{1,\Omega})={\cal U}_{n,\Omega}$ and $\mathcal{I}_i(\mathcal{H}_{1,\Omega})=\mathsf{Tr}(\mathcal{F}_2(\mathcal{H}_{1,\Omega}))=\mathsf{Tr}(\mathcal{U}_{n,\Omega})=\mathcal{U}_{1,\Omega}$. This observation proves statement 1. Analogously, to prove statement 3 notice that $\Omega$ is the only hyperedge of $\mathcal{H}_{n,\Omega}$ and thus $\mathcal{F}_i(\mathcal{H}_{n,\Omega})={\cal U}_{1,\Omega}$ and $\mathcal{I}_i(\mathcal{H}_{n,\Omega})=\mathsf{Tr}(\mathcal{U}_{1,\Omega})=\mathcal{U}_{n,\Omega}$. 

Let us demonstrate statement 2. In order to prove that $\mathcal{F}_2(\mathcal{H}_{k,\Omega})=\mathcal{U}_{n-1,\Omega}$ it is enough to show that if $B\in\mathcal{U}_{n-1,\Omega}$, then $B$ is a minimal $\mathcal{R}_2$-forcing set. 

Let $B\in\mathcal{U}_{n-1,\Omega}$ be  the set of vertices initially colored black. Without loss of generality we can take $B=\{1,\ldots,n-1\}$. In accordance with $\mathcal{R}_2$ the black subset $B'=\{1,\ldots,k-1\}$ of the hyperedge $\{1,\ldots,k-1,n\}$ forces vertex $n$ to become black, because any other hyperedge containing $B'$ has all its vertices colored black. Thus $B$ is an $\mathcal{R}_2$-forcing set.

Now let us demonstrate that if $B'\varsubsetneq B$, then $B'$ is not an $\mathcal{R}_2$-forcing set. It is enough to prove this fact in the case $B'=\{1,\ldots,n-2\}$. Consider a hyperedge $A$ such that $A\cap B'\ne\emptyset$ and $A\setminus B'\ne\emptyset$. Observe that if $k=2$, then we can assume that $A$ is either $A_1=\{1,n-1\}$ or $A_2=\{1,n\}$, while if $3\leqslant k\leqslant n-2$, then without loss of generality we can assume that $A$ is either $A_1=\{1,\ldots, k-1,n-1\}$ or $A_2=\{1,\ldots, k-1,n\}$ or $A_3=\{1,\ldots, k-2,n-1,n\}$. In any case, any black subset of $A_i$ is contained in some other hyperedge with white vertices. We conclude that $\mathcal{R}_2(B')=B'$ and so $B'$ is not an $\mathcal{R}_2$-forcing set.

Finally observe that $\mathcal{I}_2(\mathcal{H}_{k,\Omega})=\mathsf{Tr}(\mathcal{F}_2(\mathcal{H}_{k,\Omega}))=\mathsf{Tr}(\mathcal{U}_{n-1,\Omega})=\mathcal{U}_{2,\Omega}$.
\end{pf}

\begin{prop}
\label{I-Uniform}
Let $\mathcal{H}_{k,\Omega}$ be a complete hypergraph on a finite set $\Omega$, $1\leqslant k\leqslant |\Omega|$. Then $\mathcal{F}_1(\mathcal{H}_{k,\Omega})=\mathcal{U}_{|\Omega|-k+1,\Omega}$ and $\mathcal{I}_1(\mathcal{H}_{k,\Omega})=\mathcal{U}_{k,\Omega}$.
\end{prop}
\begin{pf}
Let $\Omega=\{1,\ldots,n\}$ and let us demonstrate that $\mathcal{I}_1(\mathcal{H}_{k,\Omega})=\mathcal{U}_{k,\Omega}$. It is enough to show that if $X\in\mathcal{U}_{k,\Omega}$, then $X$ is a minimal $\mathcal{R}_1$-immune set. If $k=1$ or $k=|\Omega|$ the result holds trivially (as noticed in the previous proof). Hence let us assume $2\leqslant k \leqslant |\Omega|-1$. We remind that given $X\subseteq \Omega$ and $A\in E({\cal H})$, the set $\Sigma_1(X,A)$ is defined as $\Sigma_1(X,A)=\{ A'\in N(A\setminus X) \, : \,  (A'\cap X) \setminus A \neq  \emptyset \}$. In accordance with Proposition~\ref{hyper-immune} we have to show the following two facts: (1) if $X\subseteq \Omega $ is a subset with cardinality $|X|= k$, then $|\Sigma_1(X,A)|\geqslant 1$ for all $A$ such that $|A|=k$, $A\setminus X\ne\emptyset$ and $A\cap X\ne\emptyset$; (2) if $X\subseteq \Omega $ has cardinality  $1\leqslant |X|\leqslant k-1$, then there exists $A$ such that $|A|=k$, $A\setminus X\ne\emptyset$ and $A\cap X\ne\emptyset$ and such that $|\Sigma_1(X,A)|=0$.

Let us prove fact (1). Without loss of generality we  can suppose that $X=\{1, \dots,k\}$ and we can take $A=\{1, \dots, r, a_{r+1}, \dots, a_{k}\}$ where $1\leqslant r \leqslant k-1$ and $a_{r+1}, \dots, a_{k} \in \{k+1, \dots, n\}$. Then $A'=(A\cup \{r+1\}) \setminus \{r\} \in \Sigma_1(X,A)$, because $A'\in N(A\setminus X)$ and $r+1\in(A'\cap X) \setminus A$.

To prove fact (2) we can assume $X=\{1, \dots,r\}$ where $r\leqslant k-1$. Then $\Sigma_1(X,A)=\emptyset$ if $A=X\cup \{r+1, \dots, k\}$.

To conclude the proof notice that $\mathcal{F}_1(\mathcal{H}_{k,\Omega})=\mathsf{Tr}(\mathcal{I}_1(\mathcal{H}_{k,\Omega}))=\mathsf{Tr}(\mathcal{U}_{k,\Omega})=\mathcal{U}_{n-k+1,\Omega}$.
\end{pf}

\begin{rmk}
If $k=2$ and $|\Omega|=n$, then $\mathcal{H}_{2,\Omega}$ is the complete graph $K_n$ and Propositions~\ref{uni-R2} and \ref{I-Uniform} coincide; namely, $\mathcal{F}(K_n)=\mathcal{U}_{n-1,\Omega}$ and $\mathcal{I}(K_n)=\mathcal{U}_{2,\Omega}$.
\end{rmk}

By setting $k'=|\Omega|-k+1$ we deduce from Proposition~\ref{I-Uniform} that the hypergraphs $\mathcal{H}_{k',\Omega}$ and $\mathcal{H}_{k,\Omega}$ are an $\mathcal{R}_1$-forcing realization and an $\mathcal{R}_1$-immune realization of the clutter $\mathcal{U}_{k,\Omega}$, respectively. Hence the next statement holds.

\begin{thm}
\label{s-uniformes-th-R1}
 Let $\Omega$ be a finite set.  Then, for any $1\leqslant k \leqslant |\Omega|$, the $k$-uniform clutter ${\cal U}_{k,\Omega}$ is both an $\mathcal{R}_1$-forcing clutter and an $\mathcal{R}_1$-immune clutter.
\end{thm}

Now let us see that Theorem~\ref{s-uniformes-th-R1} also holds if one considers the zero forcing rule $\mathcal{R}_2$. However, since Proposition~\ref{I-Uniform} is no longer true if the rule is $\mathcal{R}_2$, the result is not deduced so straightforwardly.
 
 \begin{thm}
\label{s-uniformes-th1}
 Let $\Omega$ be a finite set.  Then, for any $1\leqslant k \leqslant |\Omega|$, the $k$-uniform clutter ${\cal U}_{k,\Omega}$ is both an $\mathcal{R}_2$-forcing clutter and an $\mathcal{R}_2$-immune clutter.
\end{thm}

\begin{pf}
Since $\textsf{Tr}({\cal U}_{k, \Omega})={\cal U}_{|\Omega|-k+1, \Omega}$ and  $\mathsf{Tr}(\mathcal{F}_2(\mathcal{H}))=\mathcal{I}_2(\mathcal{H})$, it is enough to prove that for any $1\leqslant k \leqslant |\Omega|$, the $k$-uniform clutter ${\cal U}_{k,\Omega}$ is an $\mathcal{R}_2$-forcing clutter.

Set $\Omega=\{1,2,\ldots,n\}$, let $1\leqslant k \leqslant n$ and consider the zero forcing rule $\mathcal{R}_2$. Clearly, if $\mathcal{H}$ is the hypergraph with $E(\mathcal{H})=\{\Omega\}$,  then $\mathcal{F}(\mathcal{H})={\cal U}_{1,\Omega}$; while $\mathcal{F}(\mathcal{H})={\cal U}_{n,\Omega}$ if $E(\mathcal{H})=\{\{1\},\ldots,\{n\}\}$. So, from now on we assume that $2\leqslant k \leqslant n-1$. In this case, let  $\mathcal{H}$ be the hypergraph with vertex set $\Omega$ and set of hyperedges $E(\mathcal{H})=\{A_1,A_2, \ldots,A_m\}$, where $m=\binom{n-1}{k-1}+1$, $A_1=\Omega\setminus\{1\}$ and $A_2,\ldots,A_m$ are all the subsets of the form $\{1\}\cup A'$, where $A'\in{\cal U}_{k-1,\Omega\setminus\{1\}}$. We are going to demonstrate that $\mathcal{F}(\mathcal{H})={\cal U}_{k,\Omega}$; that is, we must prove that if $X$ is a set of vertices  of size $|X|\leqslant k$,  then $X$ is a forcing set if $|X|=k$, but it is not a forcing set if $|X|\leqslant k-1$.

Let us consider first the case $|X|=k$. If $1\in X$, then $X$ is of the form  $\{1\}\cup A'$, where $A'\in{\cal U}_{k-1,\Omega\setminus\{1\}}$. In this case, by the zero forcing rule, the set of black vertices $A'$ forces all the  white vertices of $A_1$ to be become black, and so $X$ is a forcing set. Now assume $1\not\in X$. In this case, by applying the zero forcing rule, the set of black vertices $X$ forces all the  white vertices of $A_1$. Now, applying again the zero forcing rule, we have that any subset $A'\in{\cal U}_{k-1,\Omega\setminus\{1\}}$ forces vertex $1$ to become black.

Next we consider the case $|X|=m\leqslant k-1$. Since the forcing rule is increasing it is enough to prove that if  $|X|= k-1$, then $X$ is not a forcing set.

Supose first that $|X|=k-1$ and $1\not\in X$. Without loss of generality we can assume that $X=\{2,\ldots,k\}$. We must prove that if the vertices $2,\ldots,k$ are black and the vertices $1,k+1,\ldots,n$ are white, then $\mathcal{R}^\ast(X)=X$. This is clear because $\{1,2,\ldots,k\}\in E(\mathcal{H})$ and  $\{2,\ldots,k,k+1,\ldots,n\}\in E(\mathcal{H})$ (recall that $k+1\leqslant n$).

Finally, we consider the case $|X|=k-1$ and $1\in X$. Now we can suppose that $X=\{1,2,\ldots, k-1\}$, and we must demonstrate that if the vertices $1,\ldots,k-1$ are black and the vertices $k,\ldots,n$ are white, then $\mathcal{R}^\ast(X)=X$. This fact is straightforward because  $\{1, \ldots,k-1,l\}\in E(\mathcal{H})$ for $l=k,\ldots, n$ (recall that $k\leqslant n-1$).
\end{pf}

The previous proof constructs a hypergraph $\mathcal{H}$ that is an $\mathcal{R}_2$-forcing realization of $\mathcal{U}_{k,\Omega}$. By applying the construction in the case $k'=|\Omega|-k+1$ we get a hypergraph $\mathcal{H'}$ that is an $\mathcal{R}_2$-immune realization of $\mathcal{U}_{k,\Omega}$.

Finally, let us show that if one restricts to graphs, then there are graph-realizations of $\Delta$ only for some particular values of $k$.
\begin{thm}
\label{s-uniformes-th2}
 Let $\Omega$ be a finite set  and let $1\leqslant k \leqslant |\Omega|$. Then, 
 \begin{enumerate}
	\item ${\cal U}_{k,\Omega}$ is a graph-immune clutter if and only if $k=1$ or $k=2$.
	\item ${\cal U}_{k,\Omega}$ is a graph-forcing clutter if and only if $k=|\Omega|-1$ or $k=|\Omega|$.
  \end{enumerate} 
\end{thm}

\begin{pf}
First of all observe that statement~2 follows from statement~1, because $\textsf{Tr}({\cal U}_{k, \Omega})={\cal U}_{|\Omega|-k+1, \Omega}$ and  $\mathsf{Tr}(\mathcal{I}(G))=\mathcal{F}(G)$. Hence let us prove statement~1.

Set $\Omega=\{1,2,\ldots,n\}$. It is clear that $\mathcal{I}(K_n)={\cal U}_{2,\Omega}$ and that $\mathcal{I}\left(\overline{K_n}\right)={\cal U}_{1,\Omega}$, where $K_n$ and $\overline{K_n}$ are, respectively,  the complete graph and the empty graph with vertex set $\Omega$. 

Also, it is easily checked that ${\cal U}_{n,\Omega}$ is not a graph-immune clutter if $n\geqslant 3$. Indeed, if $\mathcal{I}(G)={\cal U}_{n,\Omega}$ for some graph $G$ with vertex set $\Omega$, $|\Omega|\geqslant 3$, then  each vertex $v$ must have degree one, because $\Omega\setminus\{v\}$ is not an immune set. We conclude that $G$ is a disjoint union of copies of $K_2$, and hence  $\mathcal{I}(K_n)\ne {\cal U}_{n,\Omega}$, because, for instance, each subset of two adjacent vertices is an immune set.

Thus, the proof will be complete by showing that a contradiction is achieved if we assume that there is a graph $G$ with $V(G)= \Omega=\{1,2,\ldots,n\}$  and such that $\mathcal{I}(G)={\cal U}_{k,\Omega}$ for some $3\leqslant k\leqslant n-1$. Let $G$ be such a graph. Then, from Remark~\ref{graph-immune}, we have that ${\cal U}_{k,\Omega}$ is the collection of the inclusion-minimal elements of the family $\{X\subseteq V(G): |N(v)\cap X|\ne 1\allowbreak \textrm{ for all }v\in V(G)\setminus X\}$. Let $X\subseteq V(G)$, $|X|=k-1$, and let us denote $V(G)\setminus X$ as $\overline{X}$. Notice that $|X|, |\overline{X}|\geqslant 2$, because  $3\leqslant k\leqslant n-1$. Since $X\not\in {\cal U}_{k,\Omega}$, there is a vertex $v\in \overline{X}$ such that $|N(v)\cap X|= 1$. Next we are going to prove that if $v_X$ is the vertex belonging to $N(v)\cap X$, then $N[v]=\{v_X\}\cup \overline{X}$. Indeed, let $u\ne v$ be a vertex in $\overline{X}$ and consider $X_u=\{u\}\cup X$. Since $X_u\in {\cal U}_{k,\Omega}$ and $v\not\in X_u$, then $|N(v)\cap X_u|\ne 1$ and hence $N(v)\cap X_u=\{v_X,u\}$, because $N(v)\cap X=\{v_X\}$. We conclude that $u\in N(v)$ for any $u\in\overline{X}$, $u\ne v$, and so $N[v]=\{v_X\}\cup \overline{X}$. (See Figure~\ref{fig-1}.)

%%%%%%%%%%%%%%%%%%%%%%%%%%%%%%%%%%%%%%%%%%%%%%%%%%%%%%%%%%%%%%%%%%%%%%%%%%%%%%%%%
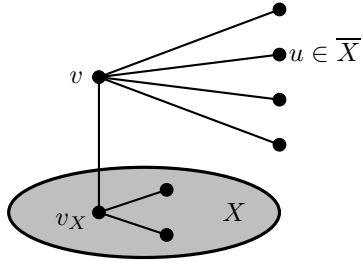
\begin{figure}[htb]
\begin{center}
\begin{tikzpicture}[scale=0.60]
%\draw [help lines] (0,0) grid (10,6);
\draw[black,thick] (4,1) ellipse (2cm and 1cm);
\filldraw[color=black, fill=gray!50, very thick] (4,1) ellipse (3cm and 1cm);
\draw[black,fill] (3,1) circle(4pt);
\draw[black,fill] (4.5,1.5) circle(4pt);
\draw[black,fill] (4.5,0.5) circle(4pt);
\draw[black,fill] (3,4) circle(4pt);
\draw[black,fill] (7,5.5) circle(4pt);
\draw[black,fill] (7,4.5) circle(4pt);
\draw[black,fill] (7,3.5) circle(4pt);
\draw[black,fill] (7,2.5) circle(4pt);
\draw[thick] (3,1)--(3,4);
\draw[thick] (3,1)--(4.5,1.5);
\draw[thick] (3,1)--(4.5,0.5);
\draw[thick] (3,4)--(7,5.5);
\draw[thick] (3,4)--(7,4.5);
\draw[thick] (3,4)--(7,3.5);
\draw[thick] (3,4)--(7,2.5);
\node at (6,1) {$X$};
\node at (2.5,4) {$v$};
\node at (2.4,0.8) {$v_X$};
\node at (8,4.6) {$u\in\overline{X}$};
\end{tikzpicture}
\end{center}
\caption{$N[v]=\{v_X\}\cup \overline{X}$}
\label{fig-1}
\end{figure}
%%%%%%%%%%%%%%%%%%%%%%%%%%%%%%%%%%%%%%%%%%%%%%%%%%%%%%%%%%%%%%%%%%%%%%%%%%%%%%%%%

Now set $X'=(X\setminus\{v_X\})\cup \{v\}$ and $\overline{X'}=V(G)\setminus X'$. As before, since $X'\not\in {\cal U}_{k,\Omega}$, there must exist a vertex $v'\in\overline{X'}$ such that $|N(v')\cap X'|= 1$. Now we are going to prove that $N[v']=N[v]$.

First let us see that $N(v')\cap X'=\{v\}$. Indeed, if $v'=v_X$, then $v'$ is not adjacent to any other vertex than $v$ of $X$ and hence $N(v')\cap X'=\{v\}$. On the other hand, if $v'\ne v_X$, then $v'\in\overline{X'}\setminus\{v_X\}= N(v)\cap\overline{X}$ and $N(v')\cap X' =\{v\}$ also holds. (See Figure~\ref{fig-2}.)

%%%%%%%%%%%%%%%%%%%%%%%%%%%%%%%%%%%%%%%%%%%%%%%%%%%%%%%%%%%%%%%%%%%%%%%%%%%%%%%%%
\begin{figure}[htb]
\begin{center}
\begin{tikzpicture}[scale=0.60]
%\draw [help lines] (0,0) grid (10,6);
%\draw[black,thick] (4,1) ellipse (2cm and 1cm);
%\draw[color=black,  very thick] (6.5,1) ellipse (3cm and 1cm);
\filldraw[rotate=-65, color=black, fill=gray!50, very thick] (-1,4.5) ellipse (3cm and 1cm);
\draw[black,fill] (8,1) circle(4pt);
\draw[black,fill] (4.5,1.5) circle(4pt);
\draw[black,fill] (4.5,0.5) circle(4pt);
\draw[black,fill] (3,4) circle(4pt);
\draw[black,fill] (7,5.5) circle(4pt);
\draw[black,fill] (7,4.5) circle(4pt);
\draw[black,fill] (7,3.5) circle(4pt);
\draw[black,fill] (7,2.5) circle(4pt);
\draw[thick] (8,1)--(3,4);
\draw[thick] (3,4)--(7,5.5);
\draw[thick] (3,4)--(7,4.5);
\draw[thick] (3,4)--(7,3.5);
\draw[thick] (3,4)--(7,2.5);
\node at (4,2.5) {$X'$};
%\node at (6.3,1) {$X$};
\node at (2.5,4) {$v$};
\node at (7.8,0.58) {$v'=v_X$};
\draw[thick] (8,1) to  (7,2.5);
\draw[thick] (8,1) to[bend right] (7,3.5);
\draw[thick] (8,1) to[bend right] (7,4.5);
\draw[thick] (8,1) to[bend right] (7,5.5);
%%%%%%%%%%%%%%%%%%%%%%%%%%%%%%%%%%
\filldraw[rotate=-65, color=black, fill=gray!50, very thick] (2.8,12.6) ellipse (3cm and 1cm);
\draw[black,fill] (16,1) circle(4pt);
\draw[black,fill] (13.5,1.5) circle(4pt);
\draw[black,fill] (13.5,0.5) circle(4pt);
\draw[black,fill] (12,4) circle(4pt);
\draw[black,fill] (16,5.5) circle(4pt);
\draw[black,fill] (16,4.5) circle(4pt);
\draw[black,fill] (16,3.5) circle(4pt);
\draw[black,fill] (17,2.5) circle(4pt);
%\draw[thick] (16,1)--(12,4);
\draw[thick] (12,4)--(16,5.5);
\draw[thick] (12,4)--(16,4.5);
\draw[thick] (12,4)--(16,3.5);
\draw[thick] (12,4)--(16,1);
\draw[thick] (12,4)--(17,2.5);
\draw[thick] (16,1)--(13.5,1.5);
\draw[thick] (16,1)--(13.5,0.5);
\node at (13,2.5) {$X'$};
\node at (11.5,4) {$v$};
\node at (16.8,0.8) {$v_X$};
\node at (17.5,2.5) {$v'$};
\draw[thick] (16,1) to[bend right] (17,2.5);
\draw[thick] (17,2.5) to[bend right] (16,3.5);
\draw[thick] (17,2.5) to[bend right] (16,4.5);
\draw[thick] (17,2.5) to[bend right] (16,5.5);
\end{tikzpicture}
\end{center}
\caption{$X'=(X\setminus\{v_X\})\cup \{v\}$, $v'=v_X$ (left), $v'\ne v_X$ (right)}
\label{fig-2}
\end{figure}
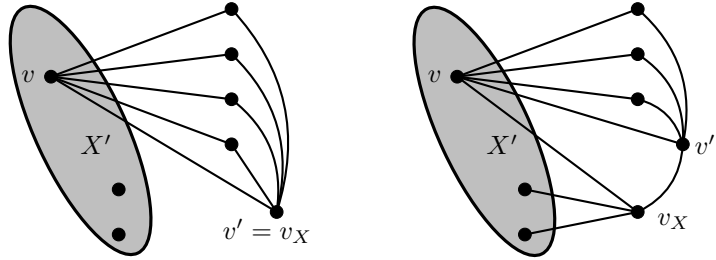
%%%%%%%%%%%%%%%%%%%%%%%%%%%%%%%%%%%%%%%%%%%%%%%%%%%%%%%%%%%%%%%%%%%%%%%%%%%%%%%%%

Now let $u\in\overline{X'}\setminus\{v'\}$ and consider $X'_u=X'\cup \{u\}$, Since $X'_u\in {\cal U}_{k,\Omega}$ and $v'\not\in X'_u$, then $|N(v')\cap X'_u|\ne 1$ and hence $N(v')\cap X'_u=\{v,u\}$, because $N(v')\cap X'=\{v\}$.  We conclude that $u\in N(v')$ for any $u\in\overline{X'}$, $u\ne v'$, and so $N[v']=\{v\}\cup \overline{X'}=\{v_X\}\cup \overline{X}=N[v]$. (See Figure~\ref{fig-2}.)

Until now we have proved that if $G$ is a graph with $V(G)=\Omega$ and such that $\mathcal{I}(G)={\cal U}_{k,\Omega}$ for some $3\leqslant k\leqslant n-1$, then there exists a pair of vertices $v$ and $v'$ such that $N[v]=N[v']$ and $|N[v]|=|N[v']|=n-k+2$. 

Let us see that the case $k=3$ is not possible. Indeed, if $\mathcal{I}(G)={\cal U}_{3,\Omega}$, then $|N[v]|=|N[v']|=n-1$, and so the only vertex $u\in\Omega\setminus N[v]=\Omega\setminus N[v']$ is not adjacent with $v$ or $v'$. Hence $N[v]=N[v']=\Omega\setminus \{u\} $. Therefore we have that $|N(u)\cap \{v,v'\}|=0$ and  $|N(u')\cap \{v,v'\}|=2$ for all $u'\in \Omega\setminus\{u,v,v'\}$. Thus we conclude that $\{v,v'\}$ is an immune set, contradicting $\mathcal{I}(G)={\cal U}_{3,\Omega}$.

At this point we can set, without loss of generality, $v=1$, $v'=2$ and $N[1]=N[2]=\{1,2,k+1,\ldots, n\}$, with $k\geqslant 4$ and $n\geqslant 5$.

Now we consider $X=\{1,2, \ldots, k\}\setminus \{3\}$. As before, since $X\not\in {\cal U}_{k,\Omega}$, there is a vertex $v\in\Omega\setminus X=\{3,k+1,\ldots,n\}$ such that $|N(v)\cap X|= 1$. Since vertices $k+1,\ldots,n$ are all adjacent to both $1$ and $2$, we conclude that $v=3$ and without loss of generality we can assume that $N(3)\cap X= \{4\}$. Clearly $N(3)\subseteq \{4,k+1, \ldots,n\}$. Let $k+1\leqslant l\leqslant n$ and consider $X_{l}=X\cup\{l\}$. Since $X_l\in {\cal U}_{k,\Omega}$ and $3\not\in X_l$, then $|N(3)\cap X_l|\ne 1$. So $N(3)\cap X_l=\{4,l\}$, because $N(3)\cap X= \{4\}$. Thus, $l\in N(3)$ if $k+1\leqslant l\leqslant n$ and hence, the equality $N[3]=\{3,4,k+1,\ldots,n\}$ holds. Furthermore, we claim that $N[4]=N[3]$. To prove our claim let us consider the set $X'=(X\cup\{3\})\setminus \{4\}$.  Since $X'\not\in {\cal U}_{k,\Omega}$, there is a vertex $v'\in\Omega\setminus X'$ such that $|N(v')\cap X'|= 1$. Now we have that $\Omega\setminus X'=\{4,k+1,\ldots,n\}$, that $N[1]=N[2]=\{1,2,k+1\ldots,n\}$,  and that $N[3]=\{3,4,k+1,\ldots,n\}$. Therefore, $v'=4$ and $N(4)\cap X'= \{3\}$ and hence  $N(4)\subseteq \{3,k+1,\ldots,n\}$. The proof of our claim will be completed by showing that $l\in N(4)$ for $k+1\leqslant l\leqslant n$. Let $k+1\leqslant l\leqslant n$ and consider $X_l'=X'\cup\{l\}$. Since $X_l'\in {\cal U}_{k,\Omega}$ and $4\not\in X_l'$, then $|N(4)\cap X_l'|\ne 1$. So $N(4)\cap X_l'=\{3,l\}$, because $N(4)\cap X'= \{3\}$. Thus, $l\in N(4)$. This completes the proof of our claim.

Finally, to achieve a contradiction, we consider the set $X=\{2,4,5,\dots, k+1\}$. Observe that $|N(v)\cap X|\ne 1$ for all $v\in \Omega \setminus X=\{1,3,k+2, \dots n\}$. We conclude from Remark~\ref{graph-immune}  that there exists $ X'\in \mathcal{I}(G)$ with $X'\subseteq X$. By assumption  $\mathcal{I}(G)={\cal U}_{k,\Omega}$. Therefore we get that $|X'|=k$ and this equality leads us to a contradiction, because $|X|=k-1$. This completes the proof of the theorem.
\end{pf}

%%%%%%%%%%%%%%%%%%%%%%%%%%%%%%%%%%%%%%%%%%%%%%%%%%%%%%%%%%%%%%%%%%%%%
\section{Conclusions} 
We have discussed two natural generalizations for hypergraphs of zero forcing on graphs such that if a graph is seen as a hypergraph, then each one of these generalized rules reduces to standard zero forcing. The families of minimal forcing and minimal immune sets of the hypergraph have been characterized within the setting of transversal theory for these forcing rules (and so also for zero forcing in graphs), and, moreover, a formulation of immune sets in terms of neighbourhoods has been provided. In Propositions~\ref{uni-R2} and \ref{I-Uniform} we have proven that the families of minimal forcing and minimal immune sets of a complete hypergraph (i.e.\  a hypergraph whose family of hyperedges is a uniform clutter) are also uniform clutters. Finally, in Theorems~\ref{s-uniformes-th-R1}, \ref{s-uniformes-th1} and \ref{s-uniformes-th2} we characterize the hypergraph and graph realizations of uniform clutters.

 A clutter $\Delta=\{A_1, \dots, A_m\}$ is $k$-\emph{homogeneous} if $|A_1|= \dots =|A_m|=k$. Hence, a clutter $\Delta$ is $k$-homogeneous if and only if $\Delta \subseteq {\cal U}_{k, \Omega}$. It is worth noting that the transversal of a homogeneous clutter is not in general homogeneous. For instance, the homogeneous clutter $\Delta=\{\{1,3\},\{1,4\},\{2,3\},\{2,4\},\{3,4\}\}$ on $\Omega=\{1,2,3,4\}$ has as transversal  the non-homogeneous clutter $\textsf{Tr}(\Delta)=\{\{1,2,3\},\{1,2,4\},\{3,4\}\}$. A continuation of the research presented in this work is the study of which homogeneous clutters have a hypergraph realization and, in such case, how many non-isomorphic realizations there are. It is worth mentioning that, since for homogeneous clutters we cannot apply transversality, the results on forcing and immune clutters will have to be handled separately.

%%%%%%%%%%%%%%%%%%%%%%%%%%%%%%%%%%%%%%%%%%%%%%%%%%%%%%%%%%%%%%%%%%%%%
\section*{Acknowledgment.} The authors would like to thank the anonymous reviewers for their careful reading and valuable suggestions.
%%%%%%%%%%%%%%%%%%%%%%%%%%%%%%%%%%%%%%%%%%%%%%%%%%%%%%%%%%%%%%%%%%%%%

%%%%%%%%%%%%%%%%%%%%%%%%%%%%%%%%%%%%%%%%%%%%%%%%%%%%%%%%%%%%%%%%%%%%%

%%%%%%%%%%%%%%%%%%%%%%%%%%%%%%%%%%%%%%%%%%%%%%%%%%%%%%%%%%%%%%%%%%%%%
%%%%%%%%%%%%%%%%%%%%%%%%%%%%%%%%%%%%%%%%%%%%%%%%%%%%%%%%%%%%%%%%%%%%%

\end{document}